\documentclass[reqno,11pt]{amsart}
\usepackage{epsfig}\usepackage{amsfonts,amsmath,amssymb,mathrsfs,verbatim,amsthm,amscd}
\newcommand{\beq}{\begin{equation}}
\newcommand{\ee}{\end{equation}}
\newcommand{\bea}{\begin{eqnarray}}
\newcommand{\eea}{\end{eqnarray}}

\usepackage{hyperref}
\def\stackreb#1#2{\ \mathrel{\mathop{#1}\limits_{#2}}}

\newcommand{\CC}{\mathbb C}
\newcommand{\R}{\mathbb R}
\newcommand{\Z}{\mathbb Z}

\newcommand{\T}{\mathbb T}

\newcommand{\ve}{\varepsilon}

\begin{document}

\vspace*{-2em}

\title[Classical special functions]
{An extended scheme of classical special functions%
}

\author{Vyacheslav \,P. Spiridonov}%

\address{Laboratory of Theoretical Physics,
Joint Institute for Nuclear Research, Dubna, Moscow region, Russia and
National Research University Higher School of Economics, Moscow, Russia
}

\maketitle

\vspace*{-3.5em}

\begin{abstract}
A unifying scheme of classical special functions of hypergeometric type obeying
orthogonality or biorthogonality relations is described. It expands the Askey scheme
of classical orthogonal polynomials and its $q$-analogue based on the Askey--Wilson polynomials.
On the top, it has two-index biorthogonal functions formed from elliptic hypergeometric
series with the absolutely continuous measure determined by the elliptic beta integral.
A new result is an inclusion of complex hypergeometric functions into the scheme.
Its further potential generalizations are discussed as well.
\end{abstract}

\vspace*{1em}

{\bf A motivation.} At some point of his career Richard Askey has proposed a catching genealogy scheme of
classical orthogonal polynomials which reached from the Chebyshev polynomials to more complicated
known families based on the simplest standard hypergeometric functions. This scheme was later coined as the
Askey scheme. About 40 years ago Andrews and Askey wrote a brief survey with a challenging title
``Classical orthogonal polynomials'' \cite{AA}. In short, they humbly claimed, by definition, the freshly
discovered Askey--Wilson polynomials \cite{AW} and all their limiting cases as ``classical'' ones.
Indeed, the Leonard characterization theorem \cite{L82} distinguished $q$-Racah polynomials
 (discrete Askey--Wilson polynomials) as the most general self-dual finite system of orthogonal polynomials.
As a result, after gathering all known limiting cases, the $q$-scheme of classical orthogonal
polynomials, or the Askey--Wilson scheme was formed with the
terminating balanced $_4\varphi_3$ basic hypergeometric series on its top \cite{KS98}.

Eight years after  \cite{AA}, Rahman and Suslov followed this trend and published a paper \cite{RS}
with a similar title ``Classical biorthogonal rational functions''. There, the set of biorthogonal
rational functions (BRF) discovered by Wilson \cite{W78,W91} and Rahman \cite{R86} (and, logically, all
their limiting cases as well) were nominated to be classical special functions.
For some time it was believed that these BRF were the most general functions of classical type,
but no rigorous justifying arguments were suggested.
And, indeed, in another eight years an elliptic extension of the Wilson BRF
was constructed by Zhedanov and the author \cite{SZ00}. A little later, an elliptic generalization of
the Rahman BRF was constructed in 2000 together with an even more general set of biorthogonal non-rational
functions \cite{Spi01,Spi03}. This note contains the author's view on the ``Grand Unification'' scheme
of classical special functions, all being related to functions of hypergeometric type,
obeying orthogonality or biorthogonality properties.
It is written by invitation of the editors of this book dedicated to R. Askey and his
influence on the development of the field of special functions.

\vspace*{1em}

{\bf Some personal reminiscences.} My initial scientific interests were associated with quantum field theory.
This had started already at my undergraduate level around 1981 from the computation in quantum chromodynamics
of the renormalization group beta function in the one-loop approximation. This continued for
about ten years further with encountering
cumbersome evaluations of many multi-loop Feynman integrals and, in time, resulted in a sharp wish
to change the subject of my scientific research. After becoming a staff member of the Baksan Neutrino Observatory
of the Institute for Nuclear Research of the Russian Academy of Sciences, I turned to theoretical aspects
of neutrino  oscillations in matter and published in 1990 a paper where
the Euler--Gauss $_2F_1$-hypergeometric function played a crucial role \cite{BS90}.
Earlier I met particular cases of this function within  educational quantum mechanics problems,
but its application to neutrino physics brought a new scientific result, which was inspiring to me.

After having worked for 30 years in the theory of special functions (accompanied by
investigations of some related problems in physics), I am computing again the endless line
of integrals. However, these investigations bring a lot more joy than my old work on Feynman integrals.
One irony is that all these pleasant integrals appeared again
to play an important role in quantum field theory --- the connection
which came to me as a complete surprise \cite{DO,Spi13}.

In 1993, when I was going from Moscow to Montr\'eal (at that time I was a postdoc at CRM),
Sergei Suslov asked me to take some of his new papers
and mail them to Askey. I did it and ``by the way'' enclosed a copy of my own paper on some
nonlinear $q$-special functions for $q^n=1,\, n\in\Z,$ appearing in quantum mechanics \cite{SS93}.
Surprisingly, I received in a short time, a letter from Askey and in it he thanked me and described
sieved polynomials (a special subclass of $q$-orthogonal polynomials for $q^n=1$), which I did not
know about \cite{AAA}. It was exactly at this time that I started to systematically learn
about orthogonal polynomials and related special functions (mostly under the influence of Alexei Zhedanov).
A friendly and equal partner
(no mentoring tone) response from the leader in the field was quite encouraging to me.

My first personal meeting with Askey happened at the first Symmetries and Integrability of
Difference Equations meeting organized by CRM near Montr\'eal in
May 1994. He firmly asked me to call him Dick and persistently stopped my attempts to use
other names in the following, especially with official titles. Ever since that time we had rare, but
regular meetings at various conferences and much more often, conversations by email.

Many details of the memory of my chats with Dick are lost in the labyrinth of my brain.
There mostly remain only certain emotional impressions from them.
At some point in a discussion with Mizan Rahman,
I expressed an admiration of Dick's enthusiasm and ability (character) to enlighten conversations beyond
purely mathematical subjects (such as education, general statements on science, Chebyshev or Ramanujan
stories, etc.). Mizan completely agreed with my feelings.

Around 1995, during the time of uncertainties, I was trying to find a position in the West.
I asked Dick to help me with that and he did not refuse.
A proposition he came up conflicted with my wish to do more science than teaching,
and I ended up declining it. Still, this was quite supportive in general and later Dick helped me with
recommendation letters for some short-term visiting positions.

The period from 2000 to 2003 was rich with personal meetings with Dick and distant communications with him.
In 2000, the Joint Institute for Nuclear Research hosted the Group Theoretical Methods in Physics Colloquium
and Dick came to Dubna for that conference as an invited speaker. This was just two months
after the very successful Special Functions-2000 two week school/workshop
at Arizona State University. So, in Dubna I had another chance to continue learning from Dick
about special functions and the people who worked on them. There is a funny
story about how he missed a connecting flight and the substitute airplane brought him
to a different Moscow airport where he was not met by the conference organizers.
After some intriguing experience of surviving in modern Russia, he managed to
reach Dubna the next day with the help of the US embassy. During the conference
we had a memorable dinner together at a cozy cafe.

In the summer of 2002, there was a large special functions meeting at the Institute for Mathematics
and Applications in Minneapolis, dedicated to the DLMF project.
The recently discovered elliptic hypergeometric functions
were at a top of interests and discussions. It was quite nice to talk about them with world experts
in the field---Andrews, Askey and many others. In particular, there I met James Wilson.
I had an interesting conversation with him on an elliptic extension of BRF which was
discovered in his PhD thesis, written under Dick's supervision.
Another memorable event happened a year later at Bexbach, Germany, where we celebrated the 70th
birthday of Dick. I sent a congratulation message to Dick on his 80th birthday and
received a warm reply from him. But unfortunately, I wasn't able to participate in the ``Askey-80'' meeting.

Another valuable story concerns the publication of the Russian edition of the book
by Andrews, Askey and Roy \cite{aar}.
In 2006, Yuri Neretin invited me to write a complementary chapter on elliptic
hypergeometric functions for the Russian translation of this book, which the authors
approved me to do.  The first draft was ready in
half a year and then there started a long battle with the publishing company on the actual publication of this
very nice textbook. I kept Dick informed on the progress and, after numerous pushes on publishers,
eventually the book appeared in print in 2013 \cite{Spi13}.
The authors were finally glad to have hard copies in their hands.

The general influence of Dick on my scientific interests is clear, but in one case it played a crucial
role. Namely, in \cite{A97} Dick described an elementary way to evaluate Rahman's $q$-beta
integral \cite{R86}, generalizing the Askey--Wilson integral \cite{AW}.
I used an elliptic extension of that method
in order to prove the elliptic beta integral evaluation formula \cite{Spi00}, which has proven to be
very helpful. At some point I said to Dick that my mathematics is actually quite elementary and
his vivid response was that he is also only able to do very simple mathematics.

\vspace*{1em}

{\bf The Askey scheme.}
Fifty years ago, the theory of orthogonal polynomials was not at the forefront of interests of the
large mathematical community. People who met orthogonal polynomials for the first time
were surprised by their amount and the list of tag-names involved in their identifications.
Therefore there was a natural hidden request to put them into some visible classification
scheme to understand relations between them. In one of his talks that I heard, Dick mentioned a
story about a conversation he once had with a person who was interested in orthogonal polynomials
(his name has dropped out of my mind). During this conversation Dick described a simple scheme
of classical orthogonal polynomials associated with the Euler--Gauss $_2F_1$ hypergeometric
function and their descendants (maybe it also included some $_3F_2$-polynomials as well,
but I do not remember this). Dick was amazed to see that this rough scheme made his companion very happy.

The simplest orthogonal polynomials are the Chebyshev polynomials of the first and second kind
which are expressible in terms of elementary functions. At the next level one has the Chebyshev--Hermite,
Charlier and Chebyshev--Laguerre polynomials expressed in terms of $_2F_0$ or $_1F_1$ functions. Going
further up in generalization, one comes to the Meixner, Krawtchouk, Legendre, ultraspherical and Meixner--Pollaczek polynomials which all use the $_2F_1$ function. Finally, the most
popular set of orthogonal polynomials carries the name of Jacobi. These polynomials are expressible
in terms of the general $_2F_1$ series and their wide usefulness is explained
by the fact that they satisfy a second order differential equation in the
argument---the hypergeometric equation.
The name combinations ``Chebyshev--Hermite'' and ``Chebyshev--Laguerre'' may sound unusual
to some people, but the point is that Chebyshev discovered in 1859 \cite{T60} these
orthogonal polynomials about four and nineteen years earlier than Hermite (1864) and Laguerre (1879),
respectively did.
Also, a more general set of discrete $_3F_2$ polynomials was described by Chebyshev already
in 1875 \cite{T75}. All these facts were mentioned by Dick many times in his talks and
to me personally, see also the explicit references in \cite{AA}.

After the discovery of Wilson \cite{W80} and, especially, Askey--Wilson \cite{AW} polynomials,
the classification scheme of orthogonal polynomials has been expanded to a
great extent. A presentation of the main part of the Askey--Wilson scheme is given in the
handbook \cite{KS98} (note that it does not include  sieved polynomials or other roots of unity systems).
Below I shall not describe its structure except for a
few elements (e.g., for brevity only Rogers polynomials are mentioned,
not the full family of $q$-analogues of the Jacobi polynomials).

The main scientific goal of this note is to sketch a more general scheme of classical
special functions. This extended scheme includes the functions that have been found to the present moment
(its further potential extensions expected in the near future are outlined in the end).
Roughly it is presented in the Table below. Let's start its
description from the top elliptic hypergeometric level which is at the right bottom of the Table.

\vskip 2mm
\centerline{\bf \scriptsize
Table 1.
\underline{\bf
CLASSICAL ORTHOGONAL AND BIORTHOGONAL SPECIAL FUNCTIONS}}
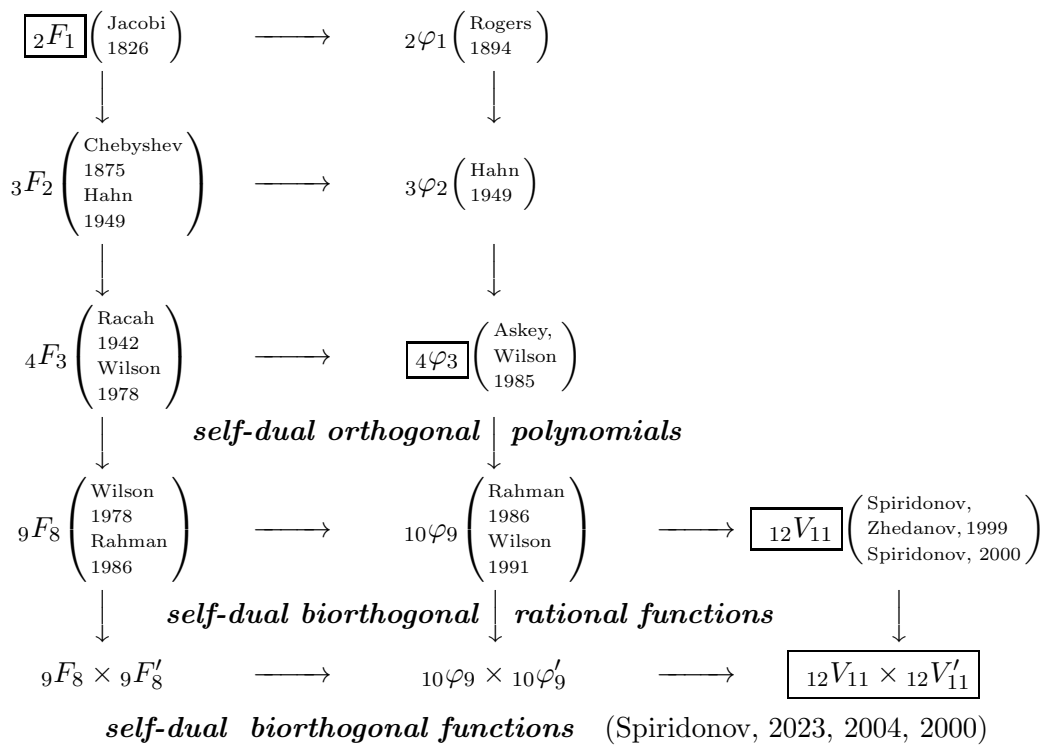
$$
\begin{CD}
\text{\hskip 5mm}
\makebox[-1.1em]{}
\fbox{{${}_2F_1$}}
{\tiny\left(\!\!
\begin{array}{l}
\text{Jacobi}\\
 \text{1826}\\
\end{array}\!\!
\right)}
@>>>
{ \makebox[-1.4em]{}\text{
${}_2\varphi_1$}}
{\tiny\left(\!\!
\begin{array}{l}
\text{Rogers} \\
\text{1894}
\end{array}\!\!
\right)}
@.\\
@VVV @VVV @.\\
{ \text{
${}_3F_2$}}
{\tiny\left(\!\!
\begin{array}{l}
\text{Chebyshev}\\
\text{1875}\\
\text{Hahn}\\
\text{1949}
\end{array}\!\!
\right)}
\text{\hskip -1mm}
@>>>
{ \text{
${}_3\varphi_2$}}
{\tiny\left(\!\!
\begin{array}{l}
\text{Hahn}\\
\text{1949}
\end{array}\!\!
\right)}
\text{\hskip 7mm}
@.\\
@VVV @VVV @.\\
{ \text{
${}_4F_3$}}
{\tiny\left(\!\!
\begin{array}{l}
\text{Racah}\\
\text{1942}\\
\text{Wilson}\\
\text{1978}
\end{array}\!\!
\right)}
@>>>
{
\fbox{${}_4\varphi_3$}}
{\tiny\left(\!\!
\begin{array}{l}
\text{Askey,}\\
\text{Wilson}\\
\text{1985}
\end{array}\!\!
\right)}
@.\\
@VVV @V{{\text{{\normalsize \bf \em  self-dual orthogonal}}\atop }\atop}
V
{{\text{{\normalsize  \bf\em polynomials}}\atop }\atop}
V@.\\
{ \text{
${}_9F_8$}}
{\tiny\left(\!\!
\begin{array}{l}
\text{Wilson}\\
\text{1978}\\
\text{Rahman}\\
\text{1986}
\end{array}\!\!
\right)}
@>>>
{ \text{
${}_{10}\varphi_9$}}
{\tiny\left(\!\!
\begin{array}{l}
\text{Rahman}\\
\text{1986}\\
\text{Wilson}\\
\text{1991}
\end{array}\!\!
\right)}
@>>>
\fbox{{ ${}_{12}V_{11}$}}
{\tiny\left(\!\!
\begin{array}{l}
\text{Spiridonov,}\\
\text{Zhedanov},
\text{1999}\\
\text{Spiridonov, 2000}
\end{array}\!\!
\right)}
@. \\
@VVV \makebox[-11em]{}  \text{\normalsize \bf \em  self-dual
biorthogonal} $\;$ @VVV
@. \makebox[-20em]{}  \text{\normalsize \bf \em \makebox[-3em]{} rational functions $\quad$ }
$\qquad$
 @   VVV \\ {}_9F_8\times {}_9F_8'
@>>>
 {}_{10}\varphi_{9}\times {}_{10}\varphi_{9}'
@>>>
\makebox[-1em]{} \fbox{{ ${}_{12}V_{11}\times {}_{12}V_{11}'$}} @.  \\
\text{\normalsize {\bf \em  \qquad\qquad self-dual }} @. \text{\normalsize {\bf \em
\makebox[-3.2em]{} biorthogonal functions \quad\;
}} @. \makebox[-7em]{} \text{\normalsize  (Spiridonov, 2023, 2004, 2000)} @.  \\
\end{CD}
$$

\vskip 5mm


{\bf The elliptic beta integral and biorthogonal functions.}
The most general special functions entering the scheme are composed
out of terminating elliptic hypergeometric series and they are biorthogonal
with respect to an absolutely continuous measure determined by the elliptic
beta integral. Therefore we start from the description of the elliptic beta
integral discovered in 2000 \cite{Spi00}.
Take 8 parameters $p, q, t_j\in\mathbb{C},\, j=1,\ldots,6,$ satisfying
the constraints $|p|, |q|, |t_j|<1$ and impose the balancing condition
$\prod_{j=1}^6t_j=pq$. Then
\begin{equation}
\frac{(p;p)_\infty(q;q)_\infty}{4\pi  i }
\int_\mathbb{T}\frac{\prod_{j=1}^6
\Gamma(t_jz^{{\pm 1}};p,q)}{\Gamma(z^{\pm 2};p,q)}\frac{dz}{z}
=\prod_{1\leq j<k\leq6}\Gamma(t_jt_k;p,q),
\label{ellbeta}\end{equation}
where $\mathbb{T}$ is the positively oriented unit circle, $(z;p)_\infty=\prod_{a=0}^\infty(1-zp^a),$ and
$$
\Gamma(z;p,q):=
\prod_{j,k=0}^\infty\frac{1-z^{-1}p^{j+1}q^{k+1}}{1-zp^{j}q^{k}}, \quad |p|, |q|<1,
$$
is the elliptic gamma function with the conventions
$$
\Gamma(t_1,\ldots,t_k;p,q):=\Gamma(t_1;p,q)\cdots\Gamma(t_k;p,q),\quad
\Gamma(tz^{\pm1};p,q):=\Gamma(tz;p,q)\Gamma(tz^{-1};p,q).
$$

The integrand function in the elliptic beta integral \eqref{ellbeta} satisfies a
linear $q$-difference equation of the first order with the coefficient given by
a particular elliptic function, which follows from the generating equation
$$
\Gamma(qz;p,q)=\theta(z;p)\Gamma(z;p,q), \quad
\theta(z;p):=(z;p)_\infty(pz^{-1};p)_\infty,
$$
where $\theta(z;p)$ is a Jacobi theta function with a specific normalization
$$
\theta(z;p)=\frac{1}{(p;p)_\infty} \sum_{k\in\mathbb{Z}} (-1)^kp^{k(k-1)/2}z^k,
\quad \theta(pz;p)=\theta(z^{-1};p)=-z^{-1}\theta(z;p).
$$

The most popular exactly computable integral having enormous number of applications is the Gaussian
integral $\int_{\R}e^{-x^2}dx=\sqrt{\pi}$. For instance, it determines the measure in the orthogonality relation
for the Chebyshev--Hermite polynomials. However, it is only one of the limiting forms of the Euler beta integral
found earlier, which determines the orthogonality measure for the Jacobi polynomials. The Askey--Wilson
$q$-beta integral \cite{AW} is the top integral of a similar use for orthogonal polynomials.
Relation \eqref{ellbeta} absorbs all these and many other similar exactly computable integrals.
It defines the normalization for the measure of the two-index biorthogonal functions given in \cite{Spi01,Spi03}.
These are currently the most general univariate special functions with classical properties.

Identity \eqref{ellbeta} serves as a germ for the whole theory of transcendental
elliptic hypergeometric functions \cite{Spi13} (including multivariable special functions
which we will not discuss here).
From the combinatorial point of view, equality \eqref{ellbeta} contains at the bottom Newton's binomial
theorem ${}_1F_0(a;x)=(1-x)^{-a}$ and its other known univariate series generalizations.
As to applications in mathematical physics,
it emerges (together with the $V$-function described below) as a solution of a particular
finite-difference equation of the second order (the elliptic hypergeometric equation)
and obeys $W(E_6)$ Weyl group of symmetries.
From the integral transformations point of view, relation \eqref{ellbeta} is a particular example
of an elliptic analogue of the standard Fourier transformation. This transformation
is related to exactly solvable models of two-dimensional statistical mechanics through
solutions of the Yang--Baxter equation. Most unexpectedly, formula \eqref{ellbeta} proves the confinement
phenomenon in a special topological sector of a four-dimensional supersymmetric
field theory through the identification of superconformal indices of dual theories.
For a brief survey of these topics, see \cite{Spi13,RR}.

The explicit form of the very-well-poised elliptic hypergeometric series in the notation
introduced in \cite{Spi03,Spi13} is:
\begin{equation}
{}_{r+1}V_{r}(t_0;t_1,\ldots,t_{r-4};q,p)
= \sum_{n=0}^\infty \frac{\theta(t_0q^{2n};p)}{\theta(t_0;p)}\prod_{m=0}^{r-4}
\frac{\theta(t_m;p;q)_n}{\theta(qt_0t_m^{-1};p;q)_n}q^n,
\label{eseries} \end{equation}
where $\theta(z;p;q)_n:=\prod_{k=0}^{n-1}\theta(zq^k;p)$ is the elliptic Pochhammer symbol.
Compact notation $\theta(x,...,y;p):=\theta(x;p)\cdots\theta(y;p)$ and
$\theta(x,...,y;p;q)_n:=\theta(x;p;q)_n\cdots\theta(y;p;q)_n$ is used below.
The word combination ``elliptic series'' in this context means that the parameters in \eqref{eseries}
satisfy the balancing condition $\prod_{k=1}^{r-4}t_k=t_0^{(r-5)/2}q^{(r-7)/2}.$
Because of that, each term of this series is an elliptic function
of all its parameters. We use multiplicative notation implying that these terms
do not change under the transformations $t_j\to p^{n_j}t_j, \, q\to p^{m}q$,
for $n_j, m\in\Z$, respecting the balancing condition for odd or even $r$.
Convergence of the infinite sum of elliptic functions \eqref{eseries} in general
(i.e., beyond the obvious termination or blow up points of the terms) is not known. However,
recently some rich set of parameter values was found when it does converge \cite{KrS23}.
In most applications it appears in the terminating form reached after imposing the
constraint $t_j=q^{-N},\, N=0,1,\ldots$, for some fixed $j$.

Let's consider the following terminating elliptic hypergeometric series
\begin{equation}
R_n(z;t_1,\ldots,t_5,t_6;q,p):={}_{12}V_{11}\Big(\frac{t_4}{t_5};\frac{q}{t_1t_5},
\frac{q}{t_2t_5},\frac{q}{t_3t_5},t_4z,\frac{t_4}{z},q^{-n},\frac{pq^n}{t_5t_6};q,p\Big),
\label{Rn}\end{equation}
where $z, t_j\in\CC,\, n\in\Z_{\geq0},$ and $\prod_{k=1}^6t_k=pq$. This function is a key building element
of the biorthogonal functions described in \cite{Spi01,Spi03} (the corresponding notation passes to ours after
the changes $t_{0},\ldots, t_4\to t_{1},\ldots, t_{5},\, A\to pq/t_6$).
Since the terminating $_{r+1}V_r$-series is an elliptic
function of all its independent continuous parameters, we thus have the following
$$
R_n(pz;\underline{t};q,p)=R_n(z;\ldots,pt_j, \ldots, p^{-1}t_k,\ldots;q,p)=R_n(z;\underline{t};q,p).
$$

Clearly, $R_n$-functions are not symmetric neither in all parameters $t_j$ nor in $p$ and $q$.
Permutation symmetry in $t_1, t_2, t_3$ is evident. Permutation of parameters $t_1$ (or $t_2$, $t_3$)
and $t_4$ results simply in the multiplication of $R_n$ by a ratio of $z$-independent
elliptic Pochhammer symbols, as it follows from the Frenkel--Turaev (elliptic Bailey)
transformation for $_{12}V_{11}$ series \cite{FT}.
$R_n$ satisfy the following three-term recurrence relation established in \cite{SZ00}
(in a different notation)
\begin{eqnarray} \label{ttr}  &&  \makebox[2em]{}
(\gamma(z)-\alpha_{n+1})B(q^n/t_5t_6)\left(R_{n+1}(z;q,p)-R_n(z;q,p)\right)
\\ \nonumber && \makebox[-2em]{}
+(\gamma(z)-\beta_{n-1})B(q^{-n})
\left(R_{n-1}(z;q,p)-R_n(z;q,p)\right)
+\rho\left(\gamma(z)-\gamma(t_4)\right)R_n(z;q,p)=0,
\nonumber\end{eqnarray}
where
\begin{eqnarray}\label{B}
&& B(x)=\frac{\theta\left(x,\frac{t_4}{t_5x},
\frac{qt_4}{t_5x},\frac{qx}{t_1t_2},\frac{qx}{t_1t_3},
\frac{qx}{t_2t_3},\frac{q\eta t_6x}{p},\frac{qt_6x}{p\eta};p
\right)}{\theta\left(\frac{t_5t_6x^2}{p},\frac{qt_5t_6x^2}{p};p\right)},
\\ \label{gamma}
&& \rho=\theta\left(\frac{qt_4t_6}{p},\frac{q}{t_1t_5},
\frac{q}{t_2t_5},\frac{q}{t_3t_5},t_4\eta,\frac{t_4}{\eta};p\right),
\\ &&
\gamma(z)=\frac{\theta(z\xi,z/\xi;p)}{\theta(z\eta,z/\eta;p)},
\qquad \alpha_n=\gamma\Big(\frac{q^n}{t_5}\Big),\qquad \beta_n=\gamma\Big(\frac{q^n}{t_6}\Big).
\end{eqnarray}
Here $\xi$ and $\eta \notin \xi^{\pm1} p^{\Z}$ are arbitrary gauge parameters
and $\gamma(z)=\gamma(pz)$ is an arbitrary elliptic
function of the second order which is referred to as the elliptic grid.
$R_n$-functions do not depend on $\xi$ nor $\eta$ since
they are eliminated from the recurrence relation by application of the addition formula for theta functions.
Taking the initial condition $R_0=1$ ($R_{-1}$ can take any finite value
since  $B(q^{-n})=0$ for $n=0$), one sees that $R_n(z)$ are rational functions of
the argument $\gamma(z)$ with poles at points $\gamma(z)=\alpha_1,\ldots,\alpha_n$.
For special discrete values of one of the parameters and the variable $z$,
the functions $R_n(z;q,p)$ reduce to elliptic analogues of Wilson's finite-dimensional
$_9F_8$ and $_{10}\varphi_9$ rational functions, which were derived in \cite{SZ00}.

Passing to the polynomials $P_n(\gamma(z))=\prod_{k=1}^n(\gamma(z)-\alpha_k) R_n(z;q,p)$,
the relation \eqref{ttr} can now be rewritten as
\begin{equation}
P_{n+1}(\gamma(z))+a_n (\gamma(z)-v_n) P_n (\gamma(z))+
u_n (\gamma(z)-\alpha_n )(\gamma(z)-\beta_{n-1})P_{n-1}(\gamma(z))=0,
\label{Pttr}\end{equation}
for some coefficients $a_n, u_n, v_n$. General systems of polynomials associated with this
type of three-term recurrence relation (called $R_{II}$-type polynomials) were studied
in detail in \cite{IM95,Z99}.

An important fact is that $R_n$-functions represent particular solutions of the elliptic
hypergeometric equation provided in \cite{SpiTh,S07} in the form
\begin{eqnarray}\nonumber &&
\frac{\prod_{j=1}^8\theta(\ve_jz;p)}{\theta(z^2,qz^2;p)} (\psi(qz)-\psi(z))
+ \frac{\prod_{j=1}^8\theta(\ve_jz^{-1};p)}{\theta(z^{-2},qz^{-2};p)}(\psi(q^{-1}z)-\psi(z))
\\ && \makebox[6em]{}
+\prod_{k=1}^6\theta\left(\frac{\ve_k \ve_8}{q};p\right)\psi(z)=0,
\label{n1'}\end{eqnarray}
where $\prod_{j=1}^8\ve_j=(pq)^2$ with $\ve_7=\ve_8/q$. Indeed, comparing this equation with
the finite-difference equation (8.7) in \cite{Spi03}, one can identify the parameters
(in our notation)
$$
\ve_j=t_j,\, j=1,\ldots,4,\quad \ve_5=\frac{pq\mu}{t_5},\quad \ve_6=\frac{qt_6}{\mu},
\quad \ve_7=\frac{t_5}{q}, \quad \ve_8=t_5,
$$
where $\psi(z)=R_n(z)$ under the quantization condition $\mu=q^n$.
The property that $R_n(z)$ satisfy second order difference equations with explicit
coefficients both in $n$ and $z$ is called `self-duality' (or `bispectrality') \cite{L82}.

Now we define the companion rational functions
\begin{equation}
T_n(z;t_1,\ldots,t_5,t_6;q,p)=R_n(z;t_1,\ldots,t_6,t_5;q,p), \quad  n\in \Z_{\geq 0},
\label{Tn}\end{equation}
which are obtained from $R_n$ by permuting the parameters $t_5$ and $t_6$.
Again, these are rational functions of the variable $\gamma(z)$.

The needed general two-index functions are defined by symmetrization in $p$ and $q$:
\begin{eqnarray} \nonumber &&
R_{nm}(z):=R_n(z;\underline{t};q,p)R_m(z;\underline{t};p,q), \quad
T_{nm}(z):=T_n(z;\underline{t};q,p)T_m(z;\underline{t};p,q),
\label{Rnm}\end{eqnarray}
where $n, m\in \Z_{\geq 0}.$
Note that the $R_m$ and $T_m$ factors are invariant under the changes $z\to qz$.
The functions $R_{nm}(z)$ and $T_{nm}(z)$ are now neither rational nor elliptic.
However, they still satisfy the same difference equations
as the corresponding separate factors due to the involved termwise $p$- or $q$-ellipticity.

The following biorthogonality relation was established by the author in \cite{Spi01,Spi03}
\begin{equation}\label{ort2}
\frac{(p;p)_\infty (q;q)_\infty}{4\pi i }
\makebox[-0.5em]{}\int_{C_{mnkl}}\makebox[-1.5em]{}
T_{nl}(z)R_{mk}(z)\frac{\prod_{j=1}^6 \Gamma(t_jz^{\pm 1};p,q)}{\Gamma(z^{\pm 2};p,q)}
\frac{d z}{z}= h_{nl}\delta_{mn}\delta_{kl}\makebox[-1em]{}
\prod_{1\leq j<s\leq 6}\Gamma(t_jt_s;p,q),
\end{equation}
where $\prod_{j=1}^6 t_j=pq$ and $h_{nl}$ are the normalization constants,
\begin{eqnarray} \nonumber
&& h_{nl}= \frac{\theta(1/t_5t_6;p)
\theta(q,qt_4/t_5,pqt_4/t_6,t_1t_2,t_1t_3,t_2t_3;p;q)_n\, q^{n}}
{\theta(q^{2n}/t_5t_6;p)
\theta(1/t_4t_5,p/t_4t_6,p/t_5t_6,t_1t_4,t_2t_4,t_3t_4;p;q)_n} \\
&&\makebox[2em]{}\times  \frac{\theta(1/t_5t_6;q)
\theta(p,pt_4/t_5,pqt_4/t_6,t_1t_2,t_1t_3,t_2t_3;q;p)_l\, p^{l}}
{\theta(p^{2l}/t_5t_6;q)
\theta(1/t_4t_5,q/t_4t_6,q/t_5t_6,t_1t_4,t_2t_4,t_3t_4;q;p)_l}.
\label{norm2}\end{eqnarray}
Here $C_{mnkl}$ denotes a positively oriented contour separating the points
$z=t_jp^aq^b,$ $t_5p^{a-k}q^{b-m}, t_6p^{a-l}q^{b-n}$ for $j=1,\ldots,4,$ and
${a,b\in\Z_{\geq0}}$ from the points with $z\to z^{-1}$ reciprocal coordinates.
Note that the subfamily of these functions determined by the restriction $t_5=t_6$
satisfies the two-index plain orthogonality relation. Since the biorthogonality measure
given by the elliptic beta integral is symmetric in $t_1,\ldots,t_6$, one can permute
all these parameters in \eqref{ort2}. It is easy to relate directly to each other
biorthogonalities obtained by permuting parameters in groups $\{t_1,\ldots, t_4\}$ and
$\{t_5, t_6\}$. However, the relations resulted from permuting parameters between
these two groups are completely non-obvious.
All the involved biorthogonal functions are self-dual in the sense mentioned above.

The two-index biorthogonality relation \eqref{ort2} represented  a principally new phenomenon
in the theory of special functions.
As shown in \cite{Ros}, it describes some general decoupling property for analytic functions
of complex variables.

Now we would like to degenerate equality \eqref{ort2} to lower levels by
taking parameters to particular limiting values. First, we set $k=l=0$, i.e., consider
the standard biorthogonality relation for rational functions $R_n$ and $T_m$. Then, after substituting into
\eqref{ort2}, $t_6=pq/A$, $A=t_1\cdots t_5$, one can take the limit $p\to 0$ with other parameters being fixed.
This yields the biorthogonality relation for Rahman's BRF \cite{R86}, which we omit here.
Other types of $p\to 0$ limits, when parameters depend on $p$ in a different way,
were considered in \cite{BR15}.

Second, we take some constant $v>0$, pass to an exponential parametrization of variables
\begin{equation}
t_j=e^{-2\pi v g_j}, \qquad z=e^{-2\pi v u},\qquad p=e^{-2\pi v\omega_2}, \qquad q=e^{-2\pi v\omega_1},
\label{newpar}\end{equation}
and take the limit $v\to 0^+$. For the theta function we have the asymptotics
$$
\theta(e^{-2\pi v u};e^{-2\pi v \omega_2})\stackreb{=}{v\to 0^+}
e^{-\frac{\pi}{6\omega_2 v}}2\sin\frac{\pi u}{\omega_2}
$$
and for the elliptic gamma function \cite{Ru97}
$$
\Gamma(e^{-2\pi v u};e^{-2\pi v\omega_1},e^{-2\pi v\omega_2})
\stackreb{=}{ v \to 0^+}
e^{-\pi\frac{2u-\omega_1-\omega_2}{12v\omega_1\omega_2}}\gamma^{(2)}(u;\omega_1,\omega_2),
$$
where $\gamma^{(2)}$ is the Faddeev's modular dilogarithm, or hyperbolic gamma function
$$
\gamma^{(2)}(u;\mathbf{\omega})= e^{-\frac{\pi i }{2}
B_{2,2}(u;\mathbf{\omega}) } \frac{(\tilde q e^{2\pi  i  \frac{u}{\omega_1}};\tilde q)_\infty}
{(e^{2\pi  i  \frac{u}{\omega_2}};q)_\infty},
\quad q=e^{2\pi i \frac{\omega_1}{\omega_2}}, \quad
\tilde q= e^{-2\pi i \frac{\omega_2}{\omega_1}},
$$
$$
 B_{2,2}(u;\mathbf{\omega})=\frac{1}{\omega_1\omega_2}
\left(\left(u-\frac{\omega_1+\omega_2}{2}\right)^2-\frac{\omega_1^2+\omega_2^2}{12}\right).
$$
Although the expression given for the $\gamma^{(2)}$-function assumes $|q|<1$, there is an
integral representation for it which is well defined for $|q|=1$ as well.

As shown in \cite{rai:limits}, both these limits are uniform on  compacta. Since in the biorthogonality relation
\eqref{ort2}, $T_{nl}(z)$ and $R_{mk}(z)$ remain bounded on the (infinite) domain of integration,
the uniformness is preserved. In this way we obtain the following two-index biorthogonality relation
(it was described in \cite{SpiTh} using an alternative approach)
\begin{eqnarray}
\int_{C_{mnkl}}T_{nl}(u)R_{mk}(u) \frac{
\prod_{j=1}^6\gamma^{(2)}(g_j\pm u;\mathbf{\omega})}{\gamma^{(2)}(\pm 2u;\mathbf{\omega})}
\frac{du}{2 i \sqrt{\omega_1\omega_2}}
= h_{nl}\delta_{mn}\delta_{kl}\makebox[-0.5em]{}
\prod_{1\leq j<k\leq 6}\gamma^{(2)}(g_j+g_k;\mathbf{\omega}),
\label{ort_hyp}\end{eqnarray}
where the balancing condition is converted to $\sum_{k=1}^6g_k=\omega_1+\omega_2$ and
the same symbols as above are used to denote different biorthogonal functions
\begin{eqnarray*} && \makebox[4em]{}
R_{nk}(u)=R_n(u;\omega_1,\omega_2)R_k(u;\omega_2,\omega_1),
\\ &&
R_n(u;\omega_1,\omega_2)={}_{10}W_9\left(\frac{s_4}{s_5};\frac{{\bf q}}{s_1s_5},\frac{{\bf q}}{s_2s_5},
\frac{{\bf q}}{s_3s_5},s_4e^{2\pi i u},{\bf q}^{-n},\frac{{\bf q}^n}{s_5s_6} ;{\bf q},{\bf q}\right),
\\  &&
T_{nk}(u)=T_n(u;\omega_1,\omega_2)T_k(u;\omega_2,\omega_1),
\quad T_n(u;\omega_1,\omega_2)=R_n(u;\omega_1,\omega_2)\big|_{g_5 \leftrightarrow g_6}
\end{eqnarray*}
where $s_j= e^{2\pi i \frac{g_j}{\omega_2}}$,  ${\bf q}=e^{2\pi i \frac{\omega_1}{\omega_2}}$  and
$$
h_{nl}=h_n(\omega_1,\omega_2)h_l(\omega_2,\omega_1), \quad
h_n(\omega_1,\omega_2)= \frac{1-\frac{1}{s_5s_6}}{1-\frac{{\bf q}^{2n}}{s_5s_6}}
\frac{({\bf q},\frac{{\bf q}s_4}{s_5},\frac{{\bf q}s_4}{s_6},s_1s_2,s_1s_3,s_2s_3;{\bf q})_n\, {\bf q}^{-n}}
{(\frac{1}{s_4s_5},\frac{1}{s_4s_6},\frac{1}{s_5s_6},s_1s_4,s_2s_4,s_3s_4;{\bf q})_n}.
$$
Here ${}_{10}W_9$ is the very-well-poised ${}_{10}\varphi_9$ $\bf q$-hypergeometric series satisfying
additionally the balancing condition.
The contour $C_{mnkl}$ goes from $- i \infty$ to $+ i \infty$ separating points
$u=g_{1,2,3,4}+\omega_2 a+\omega_1 b,$ $g_5+\omega_2(a-k)+\omega_1(b-m), g_6+\omega_2(a-l)+\omega_1(b-n)$ for
${a,b\in\Z_{\geq0}}$ from their $u\to -u$ reciprocals.

The new result of the present note consists in the reduction of the biorthogonality \eqref{ort_hyp}
to the level of complex hypergeometric functions.
In \cite{SS20b} the following asymptotic relation for the hyperbolic gamma function was rigorously
proven and shown to be uniform on compacta
\beq\label{gam2lim2}
\gamma^{(2)}( i \sqrt{\omega_1\omega_2}(n+x\delta);\omega_1,\omega_2)\stackreb{=}{\delta\to 0^+} e^{\frac{\pi i }{2}n^2} (4\pi\delta)^{ i x-1}{\bf \Gamma}(x,n),
\quad \sqrt{\omega_1\over \omega_2}= i +\delta,
\ee
where $n\in \Z, \, x\in\CC$, and  ${\bf\Gamma}(x,n)$ is the gamma function over the field of complex numbers
\begin{equation}
{\bf \Gamma}(x,n) =\frac{\Gamma(\frac{n+ i x}{2})}{\Gamma(1+\frac{n- i x}{2})}.
\label{Cgamma}\end{equation}
Heuristically, this limit was used earlier in \cite{BMS}. Complex hypergeometric functions are constructed
either as integrals over the complex plane or, in the Mellin--Barnes type form, as infinite bilateral sums of
contour integrals of combinations of gamma functions \eqref{Cgamma}, see \cite{DM20,DS2017,Is,Neretin,SS20a}
and references therein. For example, $6j$-symbols for the principal series representations of
the SL$(2,\CC)$ group were constructed in \cite{Is,DS2017} as a triple complex plane integral
equal to a single infinite sum of univariate Mellin--Barnes type integrals.

Now we apply the degeneration limit  \eqref{gam2lim2} to the hyperbolic biorthogonality relation
\eqref{ort_hyp}. For that purpose the following substitution is performed:
$$
 g_k=i\sqrt{\omega_1\omega_2}(N_k+\delta\alpha_k), \quad
 u=i\sqrt{\omega_1\omega_2}(N+\delta y), 
$$
where $ N_k, N\in \Z+\nu,\; \nu=0, \tfrac{1}{2},\; y\in\CC,$ and $\delta\to 0$.
Let's first consider the following simplest case:
$s,l,m,n=0$ and Re$(g_k)>0$. Then we can change the integration variable to $u/\sqrt{\omega_1\omega_2}$
and take the imaginary axis as the integration contour. When $\omega_2\to -\omega_1$, the parameters
$g_k/\sqrt{\omega_1\omega_2}$ become infinitesimally close to purely imaginary integer ($\nu=0$)
or half-integer ($\nu=1/2$) values. As a result, infinitely many poles from two different sets
of the integrand poles
$$
u\in \{ g_k+\omega_1\Z_{\geq0}+ \omega_2\Z_{\geq0}\}\cup \{ -g_k-\omega_1\Z_{\geq0}- \omega_2\Z_{\geq0}\}
$$
start to pinch the integration contour at the integer ($\nu=0$) or half-integer ($\nu=1/2$)
values of the variable $u/\sqrt{\omega_1\omega_2}$, and the integral diverges.
Denote now symbolically $\Delta(u)$ a combination of hyperbolic gamma functions with the set
of poles described above. Now we can write
\bea  \nonumber &&
\int_{- i \infty}^{ i \infty}\Delta(u){du\over  i \sqrt{\omega_1\omega_2}}
=\int_{-\infty}^{\infty}\Delta(i\sqrt{\omega_1\omega_2}\, x)dx
= \sum_{N\in \Z+\nu} \int_{N-1/2}^{N+1/2}\Delta(  i \sqrt{\omega_1\omega_2}\, x)dx
\\ &&  \makebox[2em]{}
=\sum_{N\in \Z+\nu} \int_{-1/2}^{1/2} \Delta( i \sqrt{\omega_1\omega_2}(N+x))dx
=\sum_{N\in \Z+\nu} \int_{-1/2\delta}^{1/2\delta}
\delta \Delta( i \sqrt{\omega_1\omega_2}(N+y\delta))dy.
\nonumber \eea
Here we rotated the integration contour clockwise by $\pi/2$ and passed to the
infinite sum of integrals over the unit intervals with the singularities
located at their middle points. The scaling of the integration variable $x\to \delta y$
was performed to avoid pinching of the contour by poles.
Now  we see that in the limit $\delta\to 0^+$ our original univariate contour integral
over the imaginary axis was converted to a bilateral sum of integrals over $\R$.
The corresponding integrands are determined  by the asymptotic expression
$\stackreb{\lim}{\delta\to 0}\,\delta\Delta( i \sqrt{\omega_1\omega_2}(N+y\delta))$,
provided it is well defined and the limit is uniform.
For sufficiently large values of indices $s,l,m,n\neq 0$ in the hyperbolic biorthogonality relation
\eqref{ort_hyp}, the integration contour is a Mellin--Barnes type deformation of the imaginary axis.
Therefore for performing the limit $\omega_2\to -\omega_1$, we should split the $u$-integration contour
into the sum of deformed finite intervals of integration attached
by the endpoints to the imaginary axis at the half-integer points. Then in the $\delta\to 0^+$ limit
these finite intervals become appropriate deformations of $\R$ for the $y$-variable integrations.

In the same way as in \cite{SS20a}, these manipulations result in diverging expressions
on both sides of \eqref{ort_hyp} and, after cancelling the common factor $(-1)^{2\nu}/(4\pi\delta)^5$,
the following biorthogonality relation is obtained (which was announced in \cite{Spi23}):
\begin{eqnarray}\nonumber &&
\frac{1}{8\pi}\sum_{N\in \Z+\nu}\int_{C_{mnsl}}T_{nl}(y,N)R_{ms}(y,N)(y^2+N^2)
\prod_{k=1}^6{\bf \Gamma}(\alpha_k\pm y,N_k\pm N)dy
\\  && \makebox[4em]{}
= h_{nl}\delta_{mn}\delta_{sl}\makebox[-0.5em]{}
\prod_{1\leq j< k \leq 6}{\bf \Gamma}( \alpha_j+\alpha_k,N_j+N_k),
\label{cbiort}\end{eqnarray}
where the balancing condition for $g_k$ was split into two restrictions
$$
\sum_{k=1}^6 \alpha_k=-2 i , \qquad \sum_{k=1}^6 N_k=0.
$$
The two-index $R_{ns}$-functions have the form
$$
R_{ns}(y,N)=R_n(Y;\underline{a})R_s(Y';\underline{a'}), \quad \underline{a}=(a_1,\ldots,a_6),
$$
$$
 Y=\frac{ i  y+N}{2}, \quad Y'=\frac{ i  y-N}{2}, \quad a_k=\frac{ i  \alpha_k+N_k}{2},
 \quad a_k'=\frac{ i  \alpha_k-N_k}{2},
$$
and they are parametrized by
the very-well-poised 2-balanced ${}_9F_8$ standard hypergeometric series
\begin{eqnarray}\nonumber && \makebox[-2em]{}
R_n(Y,\underline{a})={{}_9F_8}^{vwp}\left(a_4-a_5;1-a_1-a_5,1-a_2-a_5,1-a_3-a_5,a_4\pm Y,-n,n-a_5-a_6;1\right)
\\  \nonumber && \makebox[4em]{}
=\sum_{k=0}^{n}\frac{2k+a_4-a_5}{a_4-a_5}
\frac{(a_4-a_5,1-a_1-a_5,1-a_2-a_5,1-a_3-a_5)_k}{(1,a_1+a_4,a_2+a_4,a_3+a_4)_k}
\\ && \makebox[5em]{} \times
\frac{(a_4+Y,a_4-Y,-n,n-a_5-a_6)_k}{(1-a_5-Y,(1-a_5+Y,1+n+a_4-a_5,1-n+a_4+a_6)_k},
\label{9F8}\end{eqnarray}
where $(a_1,\ldots,a_k)_n:=(a_1)_n\cdots (a_k)_n$, $(a)_n=a(a+1)\cdots (a+n-1)$.
The companion functions are obtained by permutation of parameters
$$
T_{nl}(y,N)=T_n(Y;\underline{a})T_l(Y',\underline{a}'),
\quad T_n(Y;\underline{a})=R_n(Y;\underline{a})\big|_{N_5 \leftrightarrow N_6 \atop \alpha_5 \leftrightarrow \alpha_6},
$$
and the normalization constants read $h_{nl}=h_n(\underline{a})h_l(\underline{a}'),$ with
$$
h_n(\underline{a})= \frac{(-1)^{n}(a_5+a_6)}{a_5+a_6 -2n}
\frac{(1,1+a_4-a_5,1+a_4-a_6, a_1+a_2,a_1+a_3,a_2+a_3)_n}
{(-a_4-a_5,-a_4-a_6,-a_5-a_6, a_1+a_4,a_2+a_4,a_3+a_4)_n}.
$$
The integration contour $C_{mnsl}$ is a deformation of $\R$ which should
separate the points $y=\alpha_{1,2,3,4}-ic, \alpha_5+i(s+m-c), \alpha_6+i(l+n-c)$,
$c\in\Z_{\geq 0}$, from their $y\to -y$ reciprocals.

For $s,l,m,n=0$, equality \eqref{cbiort} is the evaluation formula for a general univariate
complex beta integral in the Mellin--Barnes representation
obtained as a limiting form of the elliptic beta integral in \cite{SS20a} and by a direct
computation in \cite{DM20}. In this case, for Im$(\alpha_k)<0$ the integration contour
can be taken as the real line $\R$.
A similar biorthogonality relation should emerge from \eqref{ort_hyp} in the
singular limit described in \cite{SS21} when the hyperbolic gamma function
degenerates to the ordinary Pochhammer symbol.

\vspace*{1em}

{\bf Classical special functions, missing items and future perspectives.}
The following elliptic analogue of the Euler--Gauss hypergeometric function
was introduced in \cite{Spi03} in the analysis of symmetries of elliptic
hypergeometric integrals
\begin{equation}
V(t_1,\dots,t_8;p,q)=\frac{(p;p)_\infty (q;q)_\infty}{4\pi i }\int_\T\frac{\prod_{j=1}^8\Gamma(t_jx^{\pm 1};p,q)}
{\Gamma(x^{\pm2};p,q)}\frac{dx}{x},
\label{V}\end{equation}
where the variables $t_j\in\CC$  are subject to the constraints $|t_j|<1$ and the balancing condition
$\prod_{j=1}^8t_j=p^2q^2$. Whenever $t_jt_k=pq$, $j\neq k,$ it reduces to the elliptic beta integral
due to the reflection equation $\Gamma(pq/z;p,q)=1/\Gamma(z;p,q)$.
This $V$-function satisfies the elliptic hypergeometric equation \eqref{n1'}
and determines its general functional solution \cite{SpiTh,S07}.
The latter  equation is the most general known finite-difference equation of the second order admitting
closed form solutions. It represents a particular contiguous relation for the $V$-function
following from the addition formula for the Jacobi theta functions and $W(E_7)$-group symmetry
transformations satisfied by \eqref{V}.
The  biorthogonal  $R_{nm}$-functions are particular limiting cases of this
function, i.e., particular solutions of the elliptic hypergeometric equation \cite{Spi03}.
Analogously, the degenerate cases of these functions satisfy corresponding degenerations
of the elliptic hypergeometric equation. At the complex hypergeometric functions level,
functions \eqref{9F8} represent particular solutions of two recurrence-difference
equations emerging in this way \cite{SS22}.
Note that the notion of self-duality gets deformed in this case. The biorthogonal functions
now satisfy second order equations in two variables (discrete and continuous ones), which are not
pure finite-difference equations in one complex indeterminate.
The existence of the general lowering operator for rational functions was analyzed in \cite{SZ07}
and the Rogrigues type relations for $R_{nm}$-functions at the elliptic level were considered by Rains
in \cite{rai:trans}. Having gathered all these data together and following the trend established in \cite{AA}
and \cite{RS}, one can call the $V$-function and all its limiting forms ``classical special functions''.
Despite Askey's persistent warnings regarding the characterization theorems as potential traps for unexpected developments, it is highly desirable to find out whether
the functions constructed in \cite{Spi01,Spi03} represent the top existing set of
special functions with the classical properties.

The scheme of classical orthogonal and biorthogonal special functions presented above is far from complete
since not all possible limiting transitions and potential generalizations have been implemented.
For instance, degenerations of trigonometric $q$-special functions
with $q=e^{2\pi i\tau}$ to the standard hypergeometric level are usually considered
in the limit $\tau\to 0$.  However, the same limit $q\to 1$ is reached at the
cusps $\tau\to \Z$ as well. In the trigonometric case it does not matter, but at the hyperbolic level,
modular transformations play a crucial role and that results in different degenerations.
When $\tau\to n,\, n\in\Z/\{0\}$, then $\tilde\tau=-1/\tau\to -1/n$
and ${\tilde q}^n=1$, and thus that particular `one half' of the hyperbolic hypergeometric integrals
degenerates to the roots of unity cases of trigonometric systems. For $n=2$ one gets
a mixture of Wilson ($q=1$) and Bannai--Ito ($\tilde q=-1$) type of orthogonal functions.
The case $n=0$ corresponds to the degeneration of hyperbolic hypergeometric functions to
the standard hypergeometric ones, since Im$(\tilde\tau)\to\infty$ (i.e., $\tilde q\to 0$).
Complex hypergeometric functions correspond to the case $n=-1$.
The singular limit $\tau\to n=1$ for the hyperbolic hypergeometric functions results in
a new class of special functions of hypergeometric type which was called rational
hypergeometric functions \cite{SS21}. An even more complicated picture can emerge in
the degenerations of hyperbolic hypergeometric functions associated with the general lens
space, since they involve a general SL$(2,\Z)$ modular transformation for $\tau$ \cite{SS20b}
mapping $\tau=n$ to arbitrary cusp.

Associated Askey--Wilson polynomials \cite{IR91,M85} are not fully considered yet.
Specifically, an important class obtained in the regime $|q|=1$ still remains unstudied.
The corresponding functions should be expressed
as a combination of two hyperbolic hypergeometric integrals of Ruijsenaars \cite{Ru99}.
For $|q|<1$ these functions are quadratic combinations of the $_8W_7$-series with the
bases $q=e^{2\pi i \tau}$ and $q=e^{-2\pi i/ \tau}$.
Similarly, the associated Rahman BRF have not been considered at all (not to mention the elliptic level).
A related problem is the description of converging infinite continued
fractions associated with the corresponding three-term recurrence relations. In particular, it is necessary
to complete the work \cite{GM98} by considering the regime $|q|=1$ when $q$ is not a root of unity.

Another remark as to incompleteness of the presented scheme is the absence of a systematic
analysis of the roots of unity cases $q^n=1$ for the trigonometric and hyperbolic hypergeometric
functions. At the Askey--Wilson orthogonal polynomials level, the first steps in that direction were
taken in \cite{AAA} and \cite{SZ97}, but the Rahman--Wilson BRF level was not touched at all.
The elliptic analogues of these roots of unity
systems are substantially more complicated being related to the commensurateness condition for
two base parameters $p^nq^m=1, \, n,m\in\Z$. They have not been properly explored yet.

A big gap lies in between the general sets of biorthogonal rational functions and orthogonal polynomials,
whose particular (classical) cases are indicated in Table 1 presented above.  Namely, the three
term recurrence relation for orthogonal polynomials differs from \eqref{Pttr} by
absence of the product $(\gamma-\alpha_n)(\gamma-\beta_{n-1})$ in front of $P_{n-1}(\gamma)$.
Applying sequential limiting procedures, it is possible to pass from this quadratic polynomial of
$\gamma$ to the linear one and further to a constant. This procedure leads to different types of
biorthogonal functions, called $R_I$-polynomials in \cite{IM95}, or Laurent biorthogonal
polynomials before reaching polynomials orthogonal on the unit circle or standard
orthogonal polynomials.

The $p\to 0$ limits for the elliptic BRF
of \cite{Spi01,Spi03} were systematically investigated in the impressive work \cite{BR15}.
A direct transition to the whole Askey--Wilson scheme
was shown to take place after different $x\to p^\alpha x$
scalings for the parameters and integration variable in the measure.
Still, the relevant considerations of that paper are nevertheless not quite satisfactory.
First of all, it is necessary to investigate
the limits of two-index biorthogonal functions which preserve this novel two-index structure, which
was skipped in \cite{BR15}. In particular, the principally important limit to the complex hypergeometric
functions preserving this two-index property in a remarkably ingenious way was missed
in \cite{rai:limits} and \cite{BR15}.
Moreover, no detailed explicit description was given for the emerging systems of special functions
lying beyond the Askey--Wilson scheme (except for the Pastro polynomials). A list of degeneration
types for the biorthogonality measure was given in \cite{BR15}. However, the explicit forms of the
corresponding biorthogonal functions themselves together with the three-term recurrence relations
and the difference equations they satisfy (like it is done in \cite{KS98}) were not presented.
Without such an analysis, it is not clear what kind of structurally different systems of functions
emerge---whether they belong to the $R_I$-family, Laurent biorthogonal polynomials, or something else.
Also, the number of derived limiting cases is  quite large and they involve complicated
technicalities. Therefore it is desirable to have an independent verification of their completeness.
Another important issue, not discussed in  \cite{BR15} at all, concerns analysis of the
degenerations for the hyperbolic BRF, which have the asymptotic behaviour qualitatively
different from the elliptic systems. As a result there should emerge hyperbolic analogues
of the limiting measures and functions considered in \cite{BR15}.
And, of course, the same questions arise with respect to the complex hypergeometric functions and standard hypergeometric functions. All this is expected to lead to an enormous, hard to grasp, enlargement of
the scheme in question.

\smallskip

In conclusion, we outline some prospects for a further hierarchic generalization of the proposed scheme.
All of these have been inspired by certain ideas from quantum field theory. Specifically, the elliptic
hypergeometric integrals on root systems have been identified with the superconformal
indices of four-dimensional minimal ${\mathcal N}=1$ supersymmetric field theories
on the compact $S^3\times S^1$ space-time (see survey \cite{RR} or the related brief description in \cite{Spi13}).
This connection provides a systematic point of view on the group-theoretical structure of most interesting
transcendental elliptic hypergeometric functions. They are defined as integrals over the Haar measure of
some compact  Lie group $G$ with the integrands given by a universal function of the
characters of a fixed set of irreducible representations of the product $G\times F$, where
$F$ is a second compact Lie group ($G$ and $F$ describe the gauge and flavour symmetries
of the field theories, respectively).

The hyperbolic hypergeometric functions were identified with the supersymmetric partition functions of
three-dimensional field theories on the squashed three-sphere $S_b^3$ with $b=\sqrt{\omega_1/\omega_2}$.
The trigonometric hypergeometric functions correspond to the space $S^2\times S^1$.
At the next stage of generalization of these special functions, one replaces $S^3$ with the general
lens space $L(r,k),\, k\in\Z_r,$ which leads to rarefied versions of hyperbolic and elliptic
hypergeometric functions. For example, the general univariate hyperbolic beta integral
corresponding to the general lens space was evaluated in \cite{SS20b}. At the elliptic level
only the simplest case $L(r,1)\times S^1$ was treated. A different situation emerges when one
computes superconformal indices of field theories over $S_b^3$, which differ from the partition functions.
The two-index biorthogonal functions constructed  by Rosengren in \cite{Ros} from
a product of two Rahman's BRF are related to the latter indices. They are expected to be connected
with the superconformal indices of four-dimensional models over $L(r,k)\times S^1$
in the limit $r\to \infty$. That is why they do not fit the scheme presented above.
The most general three- or four-dimensional manifolds for which supersymmetric field theories
can yield currently computable partition functions or superconformal indices are the
Seifert manifolds $\mathcal{M}$ or $\mathcal{M}\times S^1$, respectively.
Thus, the types of emerging special functions of hypergeometric type depend on the topology
of the manifolds over which the field theory is defined. This is a source of future generalizations
of the extended scheme of classical orthogonal functions described in this paper.

\smallskip

{\bf Acknowledgments.} I am indebted to S. E. Derkachov,
G. A. Sarkissian and A. S. Zhedanov for a long time fruitful collaboration and numerous useful discussions.
This study has been partially funded within the framework of the HSE University Basic Research Program.

\end{document}